\documentclass[12pt]{paper}

\usepackage{amssymb,amsfonts,amsthm,amsmath}

\usepackage{latexsym,multicol,graphicx}

\usepackage{graphicx, verbatim}
\usepackage{xcolor}
\usepackage{verbatim}

\usepackage[a4paper, hmargin=3cm, vmargin={4cm, 4cm}]{geometry}

\usepackage{fancyhdr}
\def\Max{\operatornamewithlimits{max\vphantom{p}}}

\newtheorem{theorem}{Theorem}
\newtheorem*{theorem*}{Theorem}
\newtheorem{lemma}[theorem]{Lemma}

\renewcommand{\phi}{\varphi}

\renewcommand{\(}{\bigl(}
\renewcommand{\)}{\bigr)\vphantom{)}}

\def\deg{\mathop{\rm deg}\nolimits}

\newcommand{\al}{\alpha}

\newcommand{\eps}{\varepsilon}

\newcommand{\la}{\lambda}

\def\done{{1\hskip-2.5pt{\rm l}}}

\renewcommand{\le}{\leqslant}
\renewcommand{\ge}{\geqslant}

\newcommand{\bR}{\mathbb R}

\newcommand{\bZ}{\mathbb Z}

\newcommand{\bT}{\mathbb T}
\newcommand{\bN}{\mathbb N}

\begin{document}

\title{Notes on the Szeg\H{o} minimum problem. \\
II. Singular measures}

\author{Alexander Borichev
\thanks{Supported by a joint grant of Russian Foundation for Basic Research
and CNRS (projects 17-51-150005-NCNI-a and PRC CNRS/RFBR 2017-2019) and by the project ANR-18-CE40-0035.} \\
\and
Anna Kononova
\thanks{Supported by a joint grant of Russian Foundation for Basic Research
and CNRS (projects 17-51-150005-NCNI-a and PRC CNRS/RFBR 2017-2019).} \\
\and
Mikhail Sodin
\thanks{Supported by ERC
Advanced Grant 692616 and ISF Grant 382/15.} }

\maketitle


\begin{abstract}
In this note, we prove several quantitative results concerning with the Szeg\H{o}
minimum problem for classes of measures on the unit circle concentrated on small subsets.
As a by-product, we refute a long-standing conjecture of Nevai.

This note can be read independently from the first one.

\end{abstract}

\section{Introduction}
In this note we will demonstrate several simple estimates of the quantity
\[
e_n(\rho)^2 = \min_{q_0, \ldots , q_{n-1}}\, \int_{\bT}
\bigl| t^n + q_{n-1}t^{n-1} + \ldots  +q_1 t + q_0 \bigr|^2\, {\rm d}\rho(t)
\]
for measures $\rho$ supported by
small subsets of the unit circle $\bT$.

We start with a straightforward lower bound for $e_n(\rho)$ for
measures $\rho$ of the form
\[
\rho = \sum_{k\ge 1} a_k \rho_k\,,
\]
where $a_k\ge 0$, $\sum_k a_k =1$, and $\rho_k$ are probability measures, $\rho_k$ is invariant w.r.t. rotation of the circle by $2\pi/2^k$ radians. This lower bound
yields a simple counter-example to the Nevai conjecture raised in~\cite{Nevai}
and then discussed by Rakhmanov in~\cite{Rakh} and by Simon
in~\cite[Sections~2.9, 9.4, 9.10]{Simon}.

Our second result (Theorem~\ref{thm_e*}) deals with discrete probability measures
\[
\rho = \sum_{j} a_j \delta_{\la_j}, \qquad \sum_j a_j =1, \ (\la_j)\subset\bT.
\]
Given a sequence $(a_j)$, we estimate the quantity
$\displaystyle \sup_{(a_j)\subset \bT} e_n(\rho) $. Our proof relies on ideas from Denisov's work~\cite{Denisov}.

Then we bring two results (Theorems~\ref{thm:superexp_metric} and~\ref{thm:superexp_capacity})
which provide conditions for super-exponential decay
of $e_n$. Note that \cite[Chapter~4]{ST}
contains a number of delicate conditions for sub-exponential decay of the sequence
$e_n(\rho)$ obtained by Erd\H{o}s--Tur\'an, Widom, Ullman, and Stahl--Totik.

We conclude this note with a discussion of the singular
continuous Riesz products for which $e_n(\rho)$
can be estimated in a simple and straightforward manner.

As in the first note, we use here the following notation: for positive $A$ and $B$, 
$A \lesssim B $ means that there is a positive numerical constant $C$ such that
$A \le CB $, while $A \gtrsim B $ means that $B \lesssim A $, and
$A \simeq B $ means that both $A \lesssim B $ and
$B \lesssim A $.

\subsubsection*{Acknowledgements}
We had several enlightening discussions of the Szeg\H{o} minimum problem with Sergei Denisov,
Leonid Golinskii, Fedor Nazarov, and Eero Saksman. It was Leonid Golinskii who told us about the Nevai conjecture. We thank all of them.

\section{Limit-invariant measures and the Nevai conjecture}

\subsection{Limit-invariant measures}
We say that a measure $\rho$ on $\bT= \bR/2\pi\bZ$ is $\alpha$-invariant if it is
invariant under the rotation $\theta\mapsto \theta + 2\pi \alpha$ $\mod 2\pi$.

\begin{lemma}\label{lemma:inv}
Let $\rho$ be a $\frac1k$-invariant measure with $k\in\bN$. Then
$e_{s}(\rho)^2 = \rho(\bT)$, $s<k$. 
\end{lemma}

\noindent{\em Proof of Lemma~\ref{lemma:inv}}:
Suppose that $k>1$ (for $k=1$ the statement is obvious). By the $\frac1k$-invariance of
the measure $\rho$, its moments of order $1\le |\ell| \le k-1$ vanish. Thus, the measures
$\rho$ and $\rho(\bT) m$ (here and elsewhere, $m$ is the normalized Lebesgue measure
on $\bT$) have the same moments of order $0\le |\ell | \le k-1$, and therefore,
\[
e_{s}(\rho)^2 = e_{s}(\rho(\bT)m)^2 = \rho(\bT) e_{s}(m)^2 = \rho(\bT),\qquad s<k,
\]
completing the proof. \hfill $\Box$

\begin{lemma}\label{lemma:limit-inv}
Suppose that $\rho$ is a probability measure on $\bT$ of the form
\[
\rho = \sum_{k\ge 1} a_k \rho_k,
\]
where $(\rho_k)$ is a sequence of probability measures such that $\rho_k$ is
$2^{-k}$-invariant, and $(a_k)$ is a sequence of non-negative numbers such
that $\sum_k a_k = 1$. Then
\[
e_{2^n} (\rho)^2 \ge \sum_{k\ge n+1} a_k.
\]
\end{lemma}

\noindent{\em Proof of Lemma~\ref{lemma:limit-inv}}: The tail
$\upsilon_n =\sum_{k\ge n+1} a_k \rho_k$ is a $2^{-(n+1)}$-invariant measure,
so that 
\[
e_{2^n}(\rho)^2 \ge e_{2^n}(\upsilon_n)^2 
=\sum_{k\ge n+1} a_k,
\]
proving the lemma. \hfill $\Box$

\medskip It is curious to observe that, generally speaking, the lower bound from
Lemma~\ref{lemma:limit-inv} cannot be significantly improved:

\begin{lemma}\label{example}
Let $\Lambda_{2^k}=\bigl\{\la\colon \la^{2^k}=1 \bigr\}$, let the sequence
$(a_k)$ be as in Lemma~\ref{lemma:limit-inv}, let
\[
\rho_k = \frac1{2^k} \sum_{\la\in\Lambda_{2^{k+1}}\setminus \Lambda_{2^k}} \delta_\la,\qquad k\ge 0,
\]
and let $\rho=\sum_{k\ge 0} a_k \rho_k$. Then
\[
\sum_{k\ge n+1} a_k \le e_{2^n}(\rho)^2 \le 4 \sum_{k\ge n} a_k\,.
\]
\end{lemma}

\noindent{\em Proof of Lemma~\ref{example}}: The measure $\rho_k$ is $2^{-k}$-invariant,
hence, the lower bound follows from Lemma~\ref{lemma:limit-inv}.

To prove the upper bound,
we put $Q_{2^n}(z)=z^{2^n}-1$. Since $Q_{2^n}$ vanishes at $\Lambda_{2^k}$ with $k\le n$
and $|Q_{2^n}|\le 2$ everywhere on $\bT$,
we have
\[
e_{2^n}^2 (\rho) \le \| Q_{2^n} \|_{L^2(\rho)}^2 \le 4 \sum_{k\ge n} a_k \rho_k (\bT)
=4 \sum_{k\ge n} a_k\,,
\]
proving the upper bound. \hfill $\Box$

\subsection{Is the relative Szeg\H{o} asymptotics always possible?}

Note that Lemma~\ref{lemma:limit-inv} yields the existence of singular
measures $\rho$ with an arbitrary slow decay of the sequence $e_n(\rho)$ (as we will
see later in Theorem~\ref{thm:RieszProd},
the Riesz products provide another construction of singular measures with such property).
Thus, taking an arbitrary measure $\mu$ with divergent logarithmic
integral
\begin{equation}\label{eq:div_integral}
\int_\bT \log \mu'\, {\rm d}m = -\infty, \quad
\mu' = \frac{{\rm d}\mu}{{\rm d}m} \,>0,
\end{equation}
and adding to $\mu$
a singular measure $\rho$ as in Lemma~\ref{lemma:limit-inv}, one can make the
sequence $e_n(\mu+\rho)$ decaying incomparably slower than the sequence $e_n(\mu)$.
It is not too difficult to achieve the same effect choosing an absolutely continuous $\rho$
such that $\mu +\rho = w\mu$ with $\log w\in L^1(m)$, or even with $\log w\in L^p(m)$ with
any $p<\infty$.

\begin{theorem}\label{thm:anti-Nevai}
Suppose that 
$\mu$ is an absolutely continuous measure on $\bT$ with
$\mu' >0$ $m$-a.e., and with divergent logarithmic integral~\eqref{eq:div_integral}.
Then, for any
sequence $\eps_n\to 0$, there exists a positive function $w$ such that, for any $p<\infty$,
$\log w\in L^p(m)$, while $e_n(w\mu)/\eps_n \to \infty$ as $n\to\infty$.
\end{theorem}

This theorem answers negatively to a question raised by Nevai in~\cite{Nevai}, where he
conjectured that {\em for any measure $\mu$ with $\mu'>0$ $m$-a.e.
and for any positive function $w$ with $\log w\in L^1(m)$, one has}
\begin{equation}
\lim_{n\to\infty} \frac{e_n(w\mu)}{e_n(\mu)}
= \exp\Bigl(\frac12 \int_\bT \log w\, {\rm d}m \Bigr).
\label{Nev}
\end{equation}
Note that when $\mu=m$ this becomes Szeg\H{o}'s theorem.
Nevai proved that this conjecture is correct when $w$ satisfies additional
regularity assumptions. Further results in that direction were obtained by
Rakhmanov~\cite{Rakh} and M\'at\'e--Nevai--Totik~\cite{MNT}. In~\cite{Rakh} (see the very
end of Section~3) Rakhmanov discusses a similar question, and guesses
that it may have a positive answer at least when $\mu$ has a smooth density
and $\log w\in L^p(m)$ with some $p>2$ (this is also refuted by Theorem~\ref{thm:anti-Nevai}). One can find a thorough discussion
of the Nevai conjecture and related topics in the Simon
treatise~\cite[Sections~2.9, 9.4, 9.10]{Simon}.

In the situation described in Theorem~\ref{thm:anti-Nevai}, relation \eqref{Nev} fails because 
for some 
unbounded $w$ with convergent logarithmic integral, we can have $e_n(w\mu)/e_n(\mu)\to\infty$, $n\to\infty$. It turns out that 
for bounded $w$ with convergent logarithmic integral and for some $\mu$, we can have $e_n(w\mu)/e_n(\mu)\to 0$, $n\to\infty$, 
which gives a different example of failure of \eqref{Nev}. 

\begin{theorem}\label{proN}
There exist an absolutely continuous measure $\mu$ and a function $w$ on $\bT$ such that $0<\mu'<1$, $0<w\le 1$ $m$-a.e., 
$\int_{\mathbb T}\log w\,\, {\rm d}m>-\infty$, and
$$
\lim_{n\to\infty}\frac{e_n(w\mu)}{e_n(\mu)}=0.
$$
\end{theorem}

\subsubsection{Proof of Theorem~\ref{thm:anti-Nevai}}
Let $\mu = e^{-H}\, m$ be a measure satisfying the
assumptions of Theorem~\ref{thm:anti-Nevai}, and set $\mu_0 = e^{-H_+}\, m\le \mu$; here and later on, $H_+=\max(H,0)$, $H_-=\max(-H,0)$. 
Then $\mu_0$ is an absolutely continuous measure on $\bT$ with
$\mu_0' >0$ $m$-a.e., and with divergent logarithmic integral~\eqref{eq:div_integral}.

The idea of the proof is straightforward: we start with the same
discrete measure $\rho$ as above, i.e.,
\[
\rho = \sum_{k\ge 1} a_k \rho_k, \qquad \rho_k =
2^{-k} \sum_{\la\in\Lambda_{2^{k+1}} \setminus \Lambda_{2^k}} \delta_\la,
\]
and spread slightly each of the measures $\rho_k$ retaining the
$2^{-k}$-invariance. First, using that $H_+<\infty$ a.e. on $\bT$, we fix $A_k$
so that
\[
m\bigl\{t\in\bT\colon |\arg (t)|<2^{-k}\pi, \
\Max_{\la\in\Lambda_{2^{k+1}} \setminus \Lambda_{2^k}} H_+(\bar\la t) > A_k \bigr\} < 2^{-k-1},
\]
and then choose a measurable set $X_k \subset \{t\in\bT\colon |\arg (t)|<2^{-k}\pi\}$
of measure $m(X_k)=\eta_k>0$ so that
\[
\sup_{t\in X_k}\, \Max_{\la\in\Lambda_{2^{k+1}} \setminus \Lambda_{2^k}} H_+(\bar\la t)
\le A_k.
\]
We choose $\eta_k$ in such a way that the sequence $(\eta_k)$ is decreasing. 

Note that, given $k$, the sets $\la X_k$, $\la\in \Lambda_{2^{k+1}}\setminus \Lambda_{2^k} $,
are disjoint. Then we set
\[
E_k = \bigcup_{\la\in\Lambda_{2^{k+1}}\setminus \Lambda_{2^k}} \la X_k,
\qquad E=\bigcup_{k\ge 1} E_k,
\]
and
\[
\widetilde\rho = \sum_{k\ge 1} a_k \widetilde{\rho}_k, \qquad \widetilde{\rho}_k =
\frac1{2^k \eta_k} \done_{E_k} \cdot m
\]
for some sequence $(a_n)$ of positive numbers to be chosen later on,  of sum $1$ 
(and observe that the measures $\widetilde\rho_k$ are $2^{-k}$-invariant probability
measures). Then we define a function $w_0$ by
\[
\mu_0 +\widetilde{\rho} = e^{-H_+} w_0\cdot m = w_0\cdot \mu_0,
\]
i.e.,
\[
w_0 = 1 + (e^{H_+} \done_E) \cdot \sum_{k\ge 1} \frac{a_k}{2^k \eta_k} \done_{E_k}.
\]

Put $w=\max(1, w_0e^{-H_-})$. Then
\[
w\cdot \mu = \max(1, w_0 e^{-H_-})e^{-H_+ + H_-} \cdot m
\ge w_0 e^{-H_+}\cdot m = w_0 \cdot \mu_0,
\]
and
\[
0 \le \log w
\le \log w_0 \le H_+ \done_E + \log_+ \Bigl( \sum_{k\ge 1} \frac{a_k}{2^k \eta_k} \done_{E_k} \Bigr)
+ \log 2.
\]

We need to choose the parameters $\eta_k$ to guarantee that both terms
on the RHS are integrable in any power $p<\infty$. Furthermore, putting
\[
\widetilde{\upsilon}_n = \sum_{k\ge n+1} a_k \widetilde{\rho}_k,
\]
recalling that the measures $\widetilde{\rho_k}$ are $2^{-k}$-invariant,
and applying Lemma~\ref{lemma:inv}, we get
\[
e_{2^n}(w\mu)^2\ge e_{2^n}(w_0\mu_0)^2 = e_{2^n} (\mu_0 +\widetilde\rho\,)^2 \ge e_{2^n}(\widetilde\rho\,)^2
\ge e_{2^n}(\widetilde\upsilon_n\,)^2 =\widetilde\upsilon_n(\mathbb T)=
\sum_{k\ge n+1} a_k.
\]

To complete the proof of Theorem~\ref{thm:anti-Nevai}, we 
choose the
sequence $a_k$ so that
\[
\eps_n = o\Bigl( \sum_{k\ge \log_2n +1} a_k \Bigr), \qquad n\to\infty\,.
\]
It remains to show that the functions $(H_+\done_E)^p$ and
$\displaystyle \log_+^p \bigl(\, \sum_{k\ge 1} \frac{a_k}{2^k \eta_k} \done_{E_k} \bigr)$
are integrable for any $p<\infty$.

We have
\[
\int_E H_+^p\, {\rm d}m  \le \sum_{k\ge 1} \int_{E_k} H_+^p \, {\rm d}m
\le \sum_{k\ge 1} A_k^p\, m(E_k) = \sum_{k\ge 1} A_k^p\, 2^k \eta_k < \infty,
\]
provided that $\eta_k$ were chosen sufficiently small with respect to $A_k$.

The second estimate is also not difficult:
\begin{align*}
\int_E \log_+^p \bigl( \sum_{k\ge 1} \frac{a_k}{2^k \eta_k} \done_{E_k} \bigr) \, {\rm d}m
&= \sum_{r\ge 1} \int_{E_r\setminus \bigcup_{s>r} E_s}
\log_+^p \bigl( \sum_{k\ge 1} \frac{a_k}{2^k \eta_k} \done_{E_k} \bigr) \, {\rm d}m \\
&= \sum_{r\ge 1} \int_{E_r\setminus \bigcup_{s>r} E_s}
\log_+^p \bigl( \sum_{k=1}^r \frac{a_k}{2^k \eta_k} \done_{E_k} \bigr) \, {\rm d}m \\
&\le \sum_{r\ge 1} \int_{E_r} \log^p \bigl( \frac{1}{\eta_r} \bigr) \, {\rm d}m \\
&\le  \sum_{r\ge 1}  2^r \eta_r\log^p \frac1{\eta_r} 
< \infty,
\end{align*}
provided that $\eta_r$ tend to zero sufficiently fast.
This finishes off the proof of Theorem~\ref{thm:anti-Nevai}. \hfill $\Box$

\subsubsection{Proof of Theorem~\ref{proN}}

Given $0<\beta<\alpha<1/2$, we set 
$h_{\alpha,\beta}(e^{2\pi {\rm i} \theta})=\alpha\done_{[0,\alpha]}(\theta)+\beta\done_{(\alpha,1/2]}(\theta)$, 
$g_{\alpha}(e^{2\pi {\rm i} \theta})=\done_{[0,\alpha]}(\theta)$.

Choose $N_k=2^{4^k}$ (so that $N_{k+1}=N_k^4$). Next, choose $\alpha_k=e^{-N_{k-2}}$, $\beta_k=e^{-N_{k+2}}$, and define 
$$
\mu=\Big(\sum_{k\ge 2}h_{\alpha_k,\beta_k}(e^{2\pi  {\rm i} N_k\theta})\Big)\cdot m.
$$

(a) Clearly, $0<\mu'<1$ $m$-a.e.\,.

(b) For every $k\ge 1$,
$$
\mu\ge \nu_k=\alpha_k\done_{[0,\alpha_k]}(e^{2\pi  {\rm i} N_k\theta})\, m.
$$
Since the measure $\nu_k$ is $1/N_k$-invariant, by Lemma~\ref{lemma:inv}, we have 
$$
e^2_s(\mu)\ge\nu_k(\mathbb T)= \alpha_k^2,\qquad 0\le s< N_k.
$$

(c) Set
$$
w(e^{2\pi  {\rm i} \theta})=\exp\Big(\sum_{k\ge 2}\log\frac{\beta_k}{\alpha_k}\cdot g_{\alpha_k}(e^{2\pi  {\rm i} N_k\theta})\Big).
$$
Then $0<w\le 1$ $m$-a.e.\ and 
$$
\int_{\mathbb T}\log(1/w)\, {\rm d}m=\sum_{k\ge 2}\alpha_k\log\frac{\alpha_k}{\beta_k}=
\sum_{k\ge 2}e^{-N_{k-2}} (N_{k+2}-N_{k-2})=\sum_{k\ge 2}N_{k-2}^{256}e^{-N_{k-2}} <\infty. 
$$

(d) Given $k\ge 3$, by construction, we have $w\mu'<2\beta_k$ on the arc 
$J=\bigl(e^{2\pi {\rm i} \theta}:1-\frac1{2N_{k-1}}<\theta<1\bigr)$ of length $\pi/N_{k-1}$ (and, in fact, on $N_{k-1}2^{-k+3}-1$ other arcs of the same length; we will not use this fact). 
Then, by \cite[Lemma 11]{BKS}, there exists a monic polynomials $T_k$ of degree $N_k$ such that
$$
|T_k(e^{2\pi {\rm i} \theta})|\le 2\cos^{N_k}\Bigl(\frac{\pi}{2N_{k-1}}\Bigr)<e^{-cN_k/N^2_{k-1}},\qquad e^{2\pi {\rm i} \theta}\in \mathbb T\setminus J.
$$
Furthermore, say, by the Remez inequality, we have 
$$
|T_k(e^{2\pi {\rm i} \theta})|\le e^{CN_k m(J)} =e^{CN_k/N_{k-1}},\qquad e^{2\pi {\rm i} \theta}\in J.
$$
Let $N_k\le n<N_{k+1}$. Then
\begin{multline*}
e_n(w\mu)\le e_{N_k}(w\mu)\le \int_{\mathbb T} |T_k|^2\, w\mu'    \, {\rm d}m
\\ \le  2\beta_k e^{CN_k/N_{k-1}}m(J)+e^{-cN_k/N^2_{k-1}}\\=
\frac{1}{N_{k-1}}e^{-N_{k+2}}e^{CN_k/N_{k-1}}+e^{-cN_k/N^2_{k-1}}
\le e^{-cN^{1/2}_k}.
\end{multline*}
On the other hand, 
$$
e_n(\mu)\ge e_{N_{k+1}-1}(\mu)\ge \alpha_{k+1}=
e^{-N_{k-1}}=e^{-N^{1/4}_k}.
$$
We conclude that 
$$
\lim_{n\to\infty}\frac{e_n(w\mu)}{e_n(\mu)}=0,
$$
which completes the proof of Theorem~\ref{proN}.
\hfill $\Box$

\section{Discrete measures on $\bT$}
Given a sequence of positive numbers
$a=(a_j)$ with $\sum_j a_j = 1$, and a sequence
$(\la_j)\subset\bT$, consider the discrete measure
\[
\rho = \sum_{j\ge 1} a_j \delta_{\la_j}.
\]
Let
\[
e^*_n(a) = \sup_{(\la_j)\subset\bT}\, e_n(\rho),
\]
and $s_k = \sum_{j>k} a_j$.

\begin{theorem}\label{thm_e*}
\mbox{}

\smallskip\noindent{\rm (i)}
Suppose that the sequence $a$ is monotonic, i.e., $a_1 \ge a_2 \ge\, \ldots\ $. Then
\[
e^*_n(a)^2 \ge (n+1)\,\sum_{j\ge 1} a_{j(n+1)}.
\]
In particular, $e_n^*(a)^2 \ge (n+1) a_{n+1}$.

\smallskip\noindent{\rm (ii)} Given $\gamma\in (0, 1)$, suppose that
\[ k |\log s_k|^{1+\frac1{\gamma}} \lesssim n.\]
Then $e_n^*(a)^2 \le C(\gamma) s_k$.

\smallskip\noindent{\rm (iii)} Given $\sigma\in (0, \tfrac12]$, suppose that
\[ k^2 | \log s_k |^{-1} \le \tfrac18 \sigma n .\]
Then $e_n^*(a)^2 \le s_k^{1-\sigma}$.
\end{theorem}
As we have already mentioned, the proofs of parts (ii) and (iii) follow ideas from Denisov's paper~\cite{Denisov}.

\subsection{Examples to Theorem~\ref{thm_e*}}
The following examples show that a combination of estimates from Theorem~\ref{thm_e*}
provides relatively tight bounds.

\subsubsection{}
Let $a=(2^{-j})_{j\ge 1}$. Then
\[
2^{-n} \le e_n^*(a)^2 \le 2^{-cn}, \quad n\in\bN.
\]

\noindent{\em Proof}: The lower bound is a straightforward consequence of (i).
To get the upper bound, we note that in this case $s_k=2^{-k}$ so we can
apply estimate (iii) with $\sigma=\tfrac12$ and $k \ge cn$.
\hfill $\Box$

\subsubsection{}
Let $a=(c(p)j^{-p})_{j\ge 1}$ with $p>1$. Then
\[
c(p) \frac1{n^{p-1}} \le e_n^*(a)^2 \le
C(p) \Bigl( \frac{\log^3 n}{n} \Bigr)^{p-1}.
\]

\noindent{\em Proof}: The lower bound is again a straightforward consequence
of (i). To prove the upper bound, first, we note that $s_k \simeq c(p) k^{1-p}$, 
so we can apply estimate (ii) with $\gamma=\tfrac12$, and $k = C(p) n (\log n)^{-3}$.
\hfill $\Box$

\medskip\noindent{\bf Remark: }
Taking $\gamma$ closer to $1$, one can improve $\log^3 n$ on the RHS
to $\log^b n$ with any $b>2$. On the other hand, it is not clear whether
the logarithmic factor is needed at all.

\subsubsection{}
Let $a=(c(p)j^{-1}\log^{-p}(j+1))_{j\ge 1}$ with
$p>1$. Then
\[
\frac{c(p)}{(\log n)^{p-1}} \le e_n^*(a)^2 \le \frac{C(p)}{(\log n)^{p-1}}.
\]

\noindent{\em Proof}: To prove the lower bound we note
that
\[
\sum_{j\ge 1} \frac1{j(n+1)\log^p(j(n+1)+1)}
\gtrsim \frac1{n \log^p n}\, \sum_{1\le j \le n} \frac1{j} \gtrsim \frac1{n\log^{p-1}n}\,.
\]
To prove the upper bound, first, we note that $s_k \ge c(p)(\log k)^{1-p}$. This
allows us to apply estimate (ii) with $\gamma=\tfrac12$,
$k = C(p) n (\log\log n)^{-3}$, for which $s_k = C(p) (\log n)^{1-p}$. \hfill $\Box$

\subsection{Proof of estimate (i)}
Consider the measure
\[
\rho =
\sum_{k=1}^{n+1}\Bigl( \sum_{j\ge 0} a_{k+j(n+1)} \Bigr) \delta_{e^{2\pi{\rm i}k/(n+1)}}\,,
\]
By the monotonicity of the sequence $a$,
\[
\min_{1\le k \le n+1}\, \sum_{j\ge 0} a_{k+j(n+1)} = \sum_{j\ge 1} a_{j(n+1)}.
\]
Hence,
\[
\rho \ge \Bigl( \sum_{j\ge 1} a_{j(n+1)} \Bigr)
\sum_{\la^{n+1}=1} \delta_\la,
\]
and Lemma~\ref{lemma:inv} yields estimate (i). \hfill $\Box$

\subsection{Proof of estimate (ii)}

Given a measure
$\rho = \sum_{j\ge 1} a_j \delta_{\la_j}$, we take $k$ and $\eps$ so that
$\eps k \ll 1 \ll \eps n$ (their values will be chosen at the end of the proof),
let $E=\{\la_1, \ldots , \la_k\}$, and, denoting by $E_{+\eps}$ the $\eps$-neighbourhood of
the set $E$, note that $m(E_{+\eps})\le 2k\eps$.

Our goal is to construct a polynomial $P$ of degree at most $n$ such that $|P(0)|\simeq 1$,
$\max_\bT |P| \lesssim 1$, and $P$ is very small on $E$. Then
\[
e_n(\rho)^2 \lesssim   \rho(\bT\setminus E)
+ \max_E |P|^2\,.
\]
The polynomial $P$ will be constructed in several steps.

\subsubsection{The outer function $F$}
Let
$ F = \exp\bigl[ -m(E_{+\eps})^{-1}
\bigl( \done_{E_{+\eps}}+{\rm i}\widetilde\done_{E_{+\eps}} \bigr)\, \bigr] $,
where $\done_{E_{+\eps}}$ is the indicator function of the set $E_{+\eps}$,
and $\widetilde\done_{E_{+\eps}}$ is its harmonic conjugate. Then, we have
\begin{itemize}
\item[(a)] $\sup_\bT |F| =1$;
\item[(b)]
$\displaystyle |F(0)| = \exp\Bigl( \int_\bT \log|F|\, {\rm d}m \Bigr) = \frac1e\,$;
\item[(c)] $ \sup_{E_{+\eps}} |F| =  \exp\bigl( -m(E_{+\eps})^{-1} \bigr)$.
\end{itemize}

\subsubsection{The trigonometric polynomial $q$ well concentrated near the origin}
Next, given $\gamma\in (0, 1)$, we
construct a trigonometric polynomial
$$q(x)=\sum_{|\ell| < n} \widehat{q}(\ell) e^{{\rm i}\ell x}$$
with the following properties:
\begin{itemize}
\item[(A)] $\widehat{q}(0) =1 $;
\item[(B)] $\displaystyle \int_{-\pi}^{\pi} |q(x)|\,{\rm d}x \le C(\gamma)$;
\item[(C)] for $s\ge 1$, $\displaystyle \int_{\frac{s}n \le |x| \le \pi} |q(x)|\,{\rm d}x \le
C(\gamma) s^{1-\gamma} e^{-s^\gamma}$.
\end{itemize}

First, we take an entire function $g$ satisfying
\[
\widehat{g}\in C_0^\infty (-1, 1), \quad
\widehat{g}(0)=1, \quad {\rm and} \quad
|g(x)|\le C(\gamma) e^{-|x|^\gamma};
\]
the construction of such entire functions is classical, see
for instance~\cite[Section~IVD]{Koosis}. Then, we let $g_n(x)=ng(nx)$, note
that the Fourier transform $\widehat{g}_n(\xi) = \widehat{g}(\xi/n)$
is supported by the interval $(-n, n)$, and consider
the periodization of $g_n$
\[
q(x) = \sum_{j\in\bZ} g_n(x-2\pi j) = \sum_{|\ell|<n} \widehat{g}(\ell/n) e^{{\rm i}\ell x}
\]
(the second equation is just the Poisson summation formula).
The RHS is a trigonometric polynomial of degree less than $n$.
It is easy to see that $q$ possesses the properties (A), (B), and (C).

\subsubsection{The algebraic polynomial $P$}
Take the Laurent polynomial $Q(e^{{\rm i}\theta}) = q(\theta)$, i.e.,
$Q(t)=\sum_{|\ell|<n} \widehat{q}(\ell) t^\ell$, and set
$P=F*Q$. This is an algebraic polynomial of degree less than $n$,
$|P(0)| = |F(0)|\cdot |\widehat{q}(0)| =e^{-1}$, and
$\max_\bT |P| \le \| F \|_{\infty, \bT} \cdot \| Q \|_{L^1(m)} \le C(\gamma)$.

To estimate $ \sup_E |P| $, we take $t=e^{{\rm i}\tau}\in E$, and proceed as follows:
\begin{align*}
|P(t)| \le \int_{-\pi}^{\pi} | F(e^{{\rm i}(\tau-\theta)}| &\cdot
|q(\theta))|\, \frac{{\rm d}\theta}{2\pi} \\
&\le \sup_\bT |F| \cdot \int_{|\theta|\ge \eps} |q|
+ \sup_{E_{+\eps}} |F| \cdot \int_{-\pi}^\pi |q| \\
&\le C(\gamma) \Bigl[ (\eps n)^{1-\gamma} e^{-(\eps n)^\gamma}
+ e^{-m(E_{+\eps})^{-1}} \Bigr].
\end{align*}
Hence,
$ \sup_E |P| \le C(\gamma) \bigl[ e^{-\frac12 (\eps n)^\gamma}
+ e^{-\frac12(\eps k)^{-1}} \bigr] $,
provided that $\eps n\ge 1$.
Thus,
\[
e_n(\rho)^2 \lesssim \max_\bT |P|^2\, \rho(\bT\setminus E)
+ \max_E |P|^2 \le C(\gamma)
\Bigl[ s_k + e^{-(\eps n)^{\gamma}} + e^{-(\eps k)^{-1}}\Bigr].
\]
At last, we set $\eps = (k |\log s_k |)^{-1}$, balancing the terms
$e^{-(\eps k)^{-1}} $ and $s_k$, and since $k |\log s_k|^{1+\frac1\gamma}
\lesssim n$, we have $e^{-(\eps n)^{\gamma}} \lesssim s_k$. \hfill $\Box$

\subsection{Proof of estimate (iii)}

Here we will use the following lemma:
\begin{lemma}[Hal\'asz~\cite{Hal}]\label{lemma:Halasz}
For any $d\in\bN$, there exists a polynomial $H_d$ of degree at most $d$ such that
$H_d(0)=1$, $H_d(1)=0$, and $\max_\bT |H_d| \le 1 +\frac2{d}$.
\end{lemma}
Note that though more general and precise estimate are known (see, for
instance,~\cite{LSV, AB}), the Hal\'asz original version suffices for
our purposes.

\medskip To prove estimate (iii), we
fix $k\le \tfrac12 n$ (to be chosen momentarily), let $d=[n/k]$, and consider the polynomial
$P(z)=\prod_{j=1}^k H_d(z\bar\la_j)$, where $H_d$ is the Hal\'asz polynomial of degree $d$
from Lemma~\ref{lemma:Halasz}. Clearly, $\deg P \le n$ and $P(0)=1$. Furthermore,
\[
\max_\bT |P| \le \Bigl( 1+ \frac2d \Bigr)^k \le e^{2k/d}
\le e^{4k^2/n}
\qquad \qquad \bigl(\, {\rm since\ } d\ge \frac{n}k-1 \ge \frac{n}{2k}\, \bigr).
\]
Thus,
\[
e_n(\rho)^2 \le \int_{\bT} |P|^2\, {\rm d}\rho \le
\bigl( \max_\bT |P|^2 \bigr) \cdot \sum_{j>k} a_j < e^{8k^2/n}\, s_k \le s_k^{1-\sigma},
\]
provided that $e^{8k^2/n} \le s_k^{-\sigma}$, that is,
$k^2/(\log s_k^{-1}) \le \tfrac18 \sigma n$. \hfill $\Box$

\section{Measures with super-exponential decay of $e_n$}
Here we bring two results, which provide conditions for super-exponential decay
of the sequence $e_n(\rho)$.

\begin{theorem}\label{thm:superexp_metric}
Let $\rho$ be a probability measure on $\bT$, red and let $n\ge 3$ be an integer.

\smallskip\noindent
{\bf (A)} Suppose that $e_n(\rho)\le e^{-\Omega}$ with $\Omega \ge 16 n\log n$.
Then there are $p\le n$ closed arcs $I_1$, \ldots , $I_p$ on $\bT$ such that
\[
\sum_{\ell =1}^p \frac{1}{\log\frac1{|I_\ell|}} \le 8\, \frac{n\log n}{\Omega}
\quad {\rm and} \quad
\rho \bigl( \mathbb T \setminus \bigcup_{1\le \ell\le p} I_{\ell} \bigr) \le e^{-\Omega}.
\]

\smallskip\noindent
{\bf (B)} Suppose that there are $p\le n/2$ closed arcs $I_1$, \ldots , $I_p$ on $\bT$
such that
\[
\sum_{\ell =1}^p \frac{1}{\log\frac1{|I_\ell|}} \le \frac{n}{2\Omega}
\quad {\rm and\quad } \rho \bigl( \mathbb T \setminus \bigcup_{1\le \ell\le p} I_\ell \bigr) \le e^{-\Omega}.
\]
Then $e_n(\rho)\le 2e^{-\frac12 \Omega}$, provided that
$\Omega\ge 4n$.
\end{theorem}

Using the logarithmic capacity (which we denote by $\rm{cap}$) we get
upper and lower bounds for $e_n(\rho)$, which are tighter than the ones
given in Theorem~\ref{thm:superexp_metric}.

\begin{theorem}\label{thm:superexp_capacity}
Let $\rho$ be a probability measure on $\bT$ and let $n\ge 2$ be a positive integer.

\smallskip\noindent{\bf (A)} Suppose that $e_n(\rho)\le e^{-\Omega}$.
Then there are $p\le n$ closed arcs $I_1$, \ldots , $I_p$ on $\bT$ such that
\[
{\rm cap} \Bigl( \bigcup_{1\le \ell \le p} I_\ell \Bigr) \le e^{-\frac12 \frac{\Omega}n}
\quad {\rm and\quad } \rho \bigl( \mathbb T \setminus \bigcup_{1\le \ell\le p} I_\ell \bigr) \le e^{-\Omega}.
\]

\smallskip\noindent
{\bf (B)} Suppose that there are $p\le n$ closed arcs $I_1$, \ldots , $I_p$
on $\bT$ such that
\[
{\rm cap} \Bigl( \bigcup_{1\le \ell \le p} I_\ell \Bigr) \le e^{-\frac{\Omega}n}
\quad
{\rm and\quad } \rho \bigl( \mathbb T \setminus \bigcup_{1\le \ell\le p} I_\ell \bigr)
\le e^{-\Omega}
\]
with $\Omega \ge C_1n$.
Then $e_{Cn}(\rho)\le e^{-\Omega/4}$.
Here $C$ and $C_1$ are positive numerical constants.
\end{theorem}

Theorem~\ref{thm:superexp_capacity} immediately yields a necessary and sufficient condition
for super-exponential decay of the sequence $e_n(\rho)$, cf.~\cite[Chapter~4]{ST}.

\begin{theorem}\label{cor-superexponent}
Let $\rho$ be a positive measure on $\bT$. Then the following
are equivalent:

\smallskip\noindent{\rm\bf (a)} the sequence $e_n(\rho)$ decays super-exponentially, i.e.,
$n^{-1}\log e_n(\rho) \to-\infty$ as $n\to\infty$;

\smallskip\noindent{\rm\bf (b)} for any positive $\eps$ and $A$, there exists $n_0$ such that
for every $n\ge n_0$ there exists a set $E\subset\bT$,  which is a union of at most $n$ arcs,
such that
\[
\operatorname{cap}(E)<\eps \quad {\rm and} \quad \rho(\bT\setminus E)<e^{-An}.
\]
\end{theorem}

\medskip\noindent{\em Proof of Theorem~\ref{cor-superexponent}:}

\smallskip\noindent{\bf (a) $\Longrightarrow$ (b)}: Suppose that the sequence $e_n(\rho)$ decays super-exponentially fast and fix $\eps$ and $A$. 
Choose $A_1\ge A$ such that $e^{-A_1/2}\le \eps$.
Then, we choose $n_0$ so that $ e_n(\rho) < e^{-A_1n} $ for $n\ge n_0$. Applying part (A) of Theorem~\ref{thm:superexp_capacity}  
with $\Omega=A_1n$, we get the set $E\subset\bT$ which
is a union of at most $n$ arcs such that $\operatorname{cap}(E)<e^{-A_1/2}\le \eps$ and
$\rho(\bT\setminus E)<e^{-A_1n}\le e^{-An}$.


\smallskip\noindent{\bf (b) $\Longrightarrow$ (a)}: Given an $A\ge C_1$ with $C_1$ as in Theorem~\ref{thm:superexp_capacity}, choose $\eps\in (0,e^{-A})$.  
By hypothesis, for every $n\ge n_0$ there exists a set $E\subset\bT$,  which is a union of at most $n$ arcs,
such that $\operatorname{cap}(E)<\eps$ and $\rho(\bT\setminus E)<e^{-An}$. Set $\Omega=An$. 
By part (B) of Theorem~\ref{thm:superexp_capacity},
for $n\ge n_0$, we have $e_{Cn}\le e^{-\Omega/4}=e^{-(A/4)n}$. 
Since $A$
can be chosen arbitrary large,
we conclude that the sequence $e_n$ decays super-exponentially fast.
\hfill $\Box$

\subsection{Proof of Theorem~\ref{thm:superexp_metric}}

\subsubsection{Proof of (A)}
Here, we will use the classical Boutroux--Cartan lower estimate of monic
polynomials outside an exceptional set.
We will bring it in the version given by Lubinsky~\cite[Theorem~2.1]{Lub}.
\begin{lemma}[Boutroux--H.~Cartan]\label{lemma:Cartan}
Given a monic polynomial $P$ of degree $n$ and an increasing sequence
$0< r_1 < r_2 < \ldots < r_n$, there exist positive integers $p\le n$ and
$(\la_j)_{j =1}^p$, $\sum_{j=1}^p \la_j =n$, and closed disks $(\bar D_j)_{j=1}^p$
of radii $2r_{\la_j}$ such that
$ \bigl\{ |P| \le \prod_{j=1}^n r_j \bigr\} \subset \bigcup_{j=1}^p \bar D_j $.
\end{lemma}
Putting $r_j=\eps j (n!)^{-1/n}$ one gets a more customary
version of this lemma~\cite[Chapter~I, Theorem~10]{Levin}, which says that
for any monic polynomial $P$ of degree $n$ and any $\eps>0$, the set
$ \bigl\{ |P| < \eps^n \bigr\}$ can be covered by at most $n$ closed disks
with the sum of radii not exceeding $2e\eps$.

\medskip
Now, turning to the proof of (A), we
suppose that $Q$ is an extremal polynomial of degree $n$. Then,
\[
e^{-2\Omega} \ge e_n^2(\rho) \ge e^{-\Omega} \rho\bigl\{ |Q|\ge e^{-\frac12 \Omega}\bigr\}\,,
\]
whence, $\rho\bigl\{ |Q|\ge e^{-\frac12 \Omega}\bigr\}\le e^{-\Omega}$.

Consider the set $\bigl\{ |Q| < e^{-\frac12 \Omega}\bigr\}$.  Put
\[
r_j = \exp\Bigl( -\frac14\, \frac{\Omega}{j \log n} \Bigr),
\quad j=1, 2, \ldots n,
\]
and note that
\[
\prod_{j=1}^n r_j = \exp\Bigl( -\frac14\, \frac{\Omega}{\log n} \sum_{j=1}^n \frac1{j} \Bigr)
> e^{-\frac12 \Omega}\,.
\]
Then, by the Bourtoux--Cartan estimate, the set
$\bigl\{ |Q| < e^{-\frac12 \Omega}\bigr\}$ can be covered by $p\le n$
arcs $I_1$, \ldots , $I_p$ of lengths $|I_\ell|=4r_{m_\ell}$, where
$\sum_\ell m_\ell = n$. Observing that
\[
4r_{m_\ell} < \exp\Bigl( -\frac14\, \frac{\Omega}{m_\ell\, \log n} + 2 \Bigr)
< \exp\Bigl( -\frac18\, \frac{\Omega}{m_\ell\, \log n} \Bigr)
\qquad ({\rm since\ } \Omega> 16 m_\ell\, \log n),
\]
we conclude that
\[
\sum_{\ell =1}^p \frac{1}{\log\frac1{|I_\ell|}} =
\sum_{\ell =1}^p \frac{1}{\log\frac1{4r_{m_\ell}}}
< \sum_{\ell=1}^p \frac{8m_\ell \log n}{\Omega} = \frac{8n\log n}{\Omega}\,,
\]
proving (A). \mbox{}\hfill $\Box$

\subsubsection{Proof of (B)}
Let $z_\ell$ be the center of the arc $I_\ell$, $\ell=1, 2, \ldots , p$.
For each $\ell$ put
\[
m_\ell = \Bigl[ \frac{\Omega}{\log\frac1{|I_\ell|}} \Bigr]
\]
and note that $ \sum_\ell m_\ell \le
\Omega\, \sum_\ell\, (\log\frac1{|I_\ell|})^{-1} \le \tfrac12 n $.
Consider the polynomial
$ P(z)=\prod_{\ell=1}^p (z-z_\ell)^{m_\ell+1} $ of degree
$ \sum_\ell m_\ell + p \le n$.
On $I_\ell$ we have
\[
|P| < 2^n |I_\ell|^{m_\ell+1} \le 2^n \exp \Bigl( \frac{\Omega}{\log\frac1{|I_\ell|}}
\cdot \log |I_\ell| \Bigr) = 2^n\, e^{-\Omega} < e^{-\frac12 \Omega}\,.
\]
Hence,
\begin{multline*}
e_n^2(\rho) \le \int_{\mathbb T} |P|^2\, {\rm d}\rho
= \Bigl( \int_{\cup_\ell I_\ell} + \int_{\mathbb T \setminus \cup_\ell I_\ell} \Bigr)
|P|^2\, {\rm d}\rho \\
\le \max_{\bigcup_\ell I_\ell} |P|^2 +
4^n \rho \bigl(\bT\setminus \bigcup_\ell I_\ell \bigr)
< e^{-\Omega} + 4^n\, e^{-\Omega}
< 2e^{-\frac12 \Omega},
\end{multline*}
proving (B). \hfill $\Box$

\subsection{Proof of Theorem~\ref{thm:superexp_capacity}}

\subsubsection{Proof of (A)}

Suppose that
$Q$ is an extremal polynomial of degree $n$ for the measure $\rho$. Then
$ \rho\bigl\{|Q| > e^{-\frac12 \Omega}\bigr\}\le e^{-\Omega}$.
Consider the set
\[
E_Q=\bigl\{|Q| \le e^{-\frac12 \Omega}\bigr\} \cap \bT
=\bigl\{|Q|^2 \le e^{-\Omega}\bigr\}\cap\bT\,.
\]
Since $|Q|^2$ is a trigonometric polynomial of degree $2n$, the set $E_Q$ is a union of
$p\le n$ closed arcs.
By a basic property of logarithmic capacity (see \cite[Theorem~5.5.4]{R}), ${\rm cap}(E_Q)\le
e^{-\frac12 \frac{\Omega}n}$.

\subsubsection{Proof of (B)}

The proof of (B) needs the following lemma.
\begin{lemma}\label{lemma:discretization}
Suppose $E\subset\bT$ is a union of at most $n\ge 14$ arcs. Then there exists
a monic polynomial $P$ of degree at most $28n$ with zeros on the unit circle such that
\[
|P|\le (\operatorname{cap}(E))^{ n}\, 2^{42n} 
\]
everywhere on $E$.
\end{lemma}

\medskip\noindent
Lemma~\ref{lemma:discretization} immediately yields (B). Indeed, for $n\ge 14$, $C=28$, and $C_1 = 80$ we have  
\begin{multline*}
 e_{28n}^2(\rho) \le \int_{\bT} |P|^2\, {\rm d}\rho
= \Bigl( \int_E + \int_{\bT\setminus E} \Bigr)  |P|^2\, {\rm d}\rho \\
\le \max_E |P|^2 + \max_\bT |P|^2\, \rho(\bT\setminus E)
\le e^{-2 \Omega}\, 4^{42n} + 4^{28n}\, e^{-\Omega} < e^{-\frac12 \Omega},
\end{multline*}
provided that $2 \cdot 4^{56n} < e^{\Omega}$. The latter condition holds whenever
{$\Omega>80 n$}. For $n<14$ we just increase $C$ and $C_1$. \hfill $\Box$

\subsubsection{Proof of Lemma~\ref{lemma:discretization}}
Let $\nu$ be the equilibrium measure of the set
$E = \bigcup_{1\le j\le p} I_j$,
$I_j=\{e^{{\rm i}\theta}\colon\alpha_j\le \theta\le \alpha'_j\}$,  $1\le j\le p\le n$, and let
\[
U^\nu (e^{{\rm i}\theta}) = \int_E \log |e^{{\rm i}\theta}-e^{{\rm i}t}|\,
{\rm d}\nu (e^{{\rm i}t})
\]
be its logarithmic potential. We assume that the measure
$\nu$ is normalized by the condition $\nu (E)=n$.
Then \[ U^\nu\bigr|_E = n \log\operatorname{cap}(E) \]
(and is $> n \log\operatorname{cap}(E)$
on $\mathbb C\setminus E$). We will construct a monic polynomial $P$ of degree
$2N$, $N \le 14n$, so that
$\log |P|\le  U^\nu + (3\log2)N$ everywhere on $E$.

For this purpose, we will replace the measure $\nu$ 
by the sum of point masses $\sum_{j=1}^{2N} \delta_{s_j}$.
It is well known (see e.g.~\cite[Lemma~4.1]{Peherstorfer--Steinbauer}
or~\cite[Lemma~3.5]{Totik}) that 
${\rm d}\nu(e^{{\rm i}\theta}) = \phi(\theta)\, {\rm d}\theta$, $e^{{\rm i}\theta} \in E$, where
\begin{equation}\label{eq:phi}
\phi(\theta) = \frac{n}{2\pi}\, \prod_{j=1}^p
\frac{|e^{{\rm i}\theta} - e^{{\rm i}\beta_j}|}{\sqrt{|e^{{\rm i}\theta}-e^{{\rm i}\alpha_j}|\cdot |e^{{\rm i}\theta} - e^{{\rm i}\alpha'_j}|}}
\end{equation}
with a sequence of points $e^{{\rm i}\beta_j}$ interlacing
with the arcs $I_j$. Since
$$
\phi(\theta)^4
= \frac{n^4}{(2\pi)^4}\prod_{j=1}^p
\frac{(e^{{\rm i}\theta}-e^{{\rm i}\beta_j})^2(1-e^{{\rm i}\theta-{\rm i}\beta_j})^2}
{(e^{{\rm i}\theta}-e^{{\rm i}\alpha_j})(1-e^{{\rm i}\theta-{\rm i}\alpha_j})(e^{{\rm i}\theta}-e^{{\rm i}\alpha'_j})(1-e^{{\rm i}\theta-{\rm i}\alpha'_j})}
$$
is a rational function of $z=e^{{\rm i}\theta}$ of degree $4p$, it has at most $8p-1$ critical points. Hence, $\phi'$ has at most $8p-1$ zeros on $[0,2\pi]$.
Thus, we can represent $E$ as a union of at most $9p-1+4n-1\le 13n-2$ arcs $\Delta'_j$,
with disjoint interiors such that $\int_{\Delta'_j}\phi\le 1/4$ and $\phi'$ has a constant sign on $\Delta'_j$. After that we split the arcs $\Delta_j'$ of length larger than or equal to $\pi/8$ 
into smaller arcs so that the length of each new arc 
is less than $\pi/8$.
Finally, we get $N\le 14n$ arcs
$\Delta_j=\{e^{{\rm i}\theta}\colon \gamma_j\le\theta\le \gamma'_j\}$ with $|\gamma_j'-\gamma_j|< \pi/8$ 
such that $\int_{\Delta_j}\phi\le 1/4$ and $\phi'$ has a constant sign on $\Delta_j$.

Set
\[
P(z)=\prod_{1\le j\le N} (z-e^{{\rm i}\gamma_j}) (z-e^{{\rm i}\gamma'_j})\,,
\qquad \deg P = 2N \le 28 n.
\]
We need to show that
\begin{equation}\label{eq:*}
\log |P(z)| \le U^\nu(z) + (3\log 2) N,  \quad z\in E\,.
\end{equation}
Fix
a point $z=e^{{\rm i}\theta}\in\Delta_j$ at which we will check this bound. Then
\begin{multline*}
\log|P(z)| = \log\bigl( |z-e^{{\rm i}\gamma_j}|\cdot |z-e^{{\rm i}\gamma_j'}| \bigr) \\
+ \Bigl(\ \sum_{\substack{{\rm dist}(z, \Delta_k) \le \frac12, \\ k\ne j}}
+  \sum_{{\rm dist}(z, \Delta_k) > \frac12} \ \Bigr)
\log\bigl( |z-e^{{\rm i}\gamma_k}|\cdot |z-e^{{\rm i}\gamma_k'}| \bigr)
\end{multline*}
The last sum does not exceed $ (\log 4) N$.

If ${\rm dist}(z, \Delta_k) \le \tfrac12$, $k\ne j$,
then $\Delta_k \subset D(z, 1)$, and
$ {\rm dist}(z, \Delta_k) = |z - e^{{\rm i}\widetilde{\gamma}_k}|$, where $\widetilde{\gamma}_k$ is one of two points $\gamma_k$, $\gamma_k'$.
Then, recalling that $\nu(\Delta_k)\le 1/4<1$ and using monotonicity of the logarithm function,
we see that
\begin{multline*}
\log\bigl( |z-e^{{\rm i}\gamma_k}|\cdot |z-e^{{\rm i}\gamma_k'}| \bigr)
\le \log |z-e^{{\rm i}\widetilde{\gamma}_k}| \\
\le \int_{\Delta_k} \log|z-e^{{\rm i}\widetilde{\gamma}_k}|\, {\rm d}\nu(e^{{\rm i}t})
\le \int_{\Delta_k} \log|z-e^{{\rm i}t}|\, {\rm d}\nu(e^{{\rm i}t}),
\end{multline*}
Hence, letting $E_0=\cup_{{\rm dist}(z, \Delta_k) \le \frac12, \, k\ne j}\Delta_k$, $E_1=E\setminus(E_0\cup \Delta_j)$, we obtain that  
\begin{align*}
\sum_{{\rm dist}(z, \Delta_k) \le \frac12, \, k\ne j}&
\log\bigl( |z-e^{{\rm i}\gamma_k}|\cdot |z-e^{{\rm i}\gamma_k'}| \bigr)\\
&\le \int_{E_0}\log |z-e^{{\rm i}t}|\, {\rm d}\nu(e^{{\rm i}t}) \\
&\le \int_{E_0}\log |z-e^{{\rm i}t}|\, {\rm d}\nu(e^{{\rm i}t}) 
\\&+
\int_{E_1\cap D(z, 1)}\log |z-e^{{\rm i}t}|\, {\rm d}\nu(e^{{\rm i}t}) + (\log 2)\nu(E_1\cap D(z, 1))\\
&+ \int_{E_1\setminus D(z, 1)}\log |z-e^{{\rm i}t}|\, {\rm d}\nu(e^{{\rm i}t}) 
\\
&\le \int_{E\setminus\Delta_j}
\log |z-e^{{\rm i}t}|\, {\rm d}\nu(e^{{\rm i}t}) + (\log 2)N.
\end{align*}

That is,
\[
\log |P(z)| \le \log\bigl( |z-e^{{\rm i}\gamma_j}|\cdot |z-e^{{\rm i}\gamma_j'}| \bigr)
+ \int_{E\setminus\Delta_j}
\log |z-e^{{\rm i}t}|\, {\rm d}\nu(e^{{\rm i}t}) + (3\log 2)N.
\]
To complete the proof of~\eqref{eq:*}, it remains to show that
\begin{equation}\label{eq:2}
\log\bigl( |z-e^{{\rm i}\gamma_j}|\cdot |z-e^{{\rm i}\gamma_j'}| \bigr)
<  \int_{\Delta_j}
\log |z-e^{{\rm i}t}|\, {\rm d}\nu(e^{{\rm i}t})\,.
\end{equation}

To do this, we are going to prove that
\begin{equation}\label{eq:3}
4\int_{\gamma_j}^{\gamma_j'} \phi (t)\, \log\frac1{|t-\theta|}\,  {\rm d}t
\le 3 \log\frac1{|\theta-\gamma_j|\cdot |\theta-\gamma_j'|} + 2,\qquad \gamma_j<\theta<\gamma_j',
\end{equation}
with the function $\phi$ defined in~\eqref{eq:phi}.

First, we verify that~\eqref{eq:3} yields~\eqref{eq:2}.
Since ${\rm d}\nu (e^{{\rm i}t})=\phi(t)\, {\rm d}t$, $\nu(\Delta_j)\le 1/4$,  and
$\tfrac1{\pi}|\theta-t| \le |e^{{\rm i}\theta}-e^{{\rm i}t}| $, $\theta,t\in \Delta_j$, 
estimate~\eqref{eq:3} yields
\[
4\int_{\Delta_j}
\log\frac1{|z-e^{{\rm i}t}|}\, {\rm d}\nu(e^{{\rm i}t})
\le 3 \log\frac1{|\theta-\gamma_j|\cdot |\theta-\gamma_j'|} + \log\pi + 2\,,
\]
where $z=e^{{\rm i}\theta}$. 
Furthermore, since the length of each arc $\Delta_j$ does not exceed $\pi/8$, we have
\[
|\theta-\gamma_j|\cdot |\theta-\gamma_j'| \le \frac14\, (\gamma_j-\gamma_j')^2
\le \frac14\,\Bigl( \frac{\pi}8\Bigr)^2,
\]
and then,
\[
\log\frac1{|\theta-\gamma_j|\cdot |\theta-\gamma_j'|}
\ge\log\Bigl( 4\cdot\Bigl( \frac8{\pi} \Bigr)^2\Bigr).
\]
Since $e^2 \cdot \pi^3 < 256$, the RHS of the last displayed formula
is bigger than $\log\pi + 2$, which gives us
\[
4\int_{\Delta_j}
\log\frac1{|z-e^{{\rm i}t}|}\, {\rm d}\nu(e^{{\rm i}t})
< 4 \log\frac1{|\theta-\gamma_j|\cdot |\theta-\gamma_j'|}
< 4 \log\frac1{|z-e^{{\rm i}\gamma_j}|\cdot |z-e^{{\rm i}\gamma_j'}|}\,,
\]
which is \eqref{eq:2}.
Thus, it remains to verify~\eqref{eq:3}.

Set $\beta=\theta-\gamma_j$, $\beta'=\gamma'_j-\theta$,
and $\psi(t)=4\phi(t+\theta)$. Then $\beta,\beta'\in(0,1)$ and 
$\displaystyle \int_{-\beta}^{\beta'} \psi (t)\, {\rm d}t\le 1$. 
We need to show that
\[
\int_{-\beta}^{\beta'} \psi (t)\, \log\frac1{|t|}\,  {\rm d}t
\le 3\log\frac1{\beta\cdot \beta'}+2\,.
\]
We assume that $\psi$ increases on $(-\beta, \beta')$, and set
$\displaystyle \psi_1(x) =\int_0^x \psi (t)\, {\rm d}t$. 
Note that the function $\psi_1$
is convex, vanishes at the origin, and $\psi_1(\beta')\le 1$, so
$ 0 \le \psi_1 (x) \le x/\beta'$ on $[0, \beta']$ and
$\psi(0)=\psi_1'(0)\le 1/\beta'$. Then, integrating by parts,
we get
\[
\int_0^{\beta'} \psi (t)\log\frac1t\, {\rm d}t
= \psi_1(\beta')\log\frac1{\beta'} + \int_0^{\beta'} \frac{\psi_1(t)}t\, {\rm d}t
\le \log\frac1{\beta'} + 1.
\]
If $\beta'<\beta$, then
\begin{multline*}
\int_0^\beta \psi(-t)\log\frac1t\, {\rm d}t
=\Bigl( \int_0^{\beta'} + \int_{\beta'}^\beta \Bigr) \psi(-t)\log\frac1t\, {\rm d}t \\
\le \psi(0) \bigl( \beta' + \beta'\log\frac1{\beta'} \bigr) + \log\frac1{\beta'}
< 2\log\frac1{\beta'} + 1,
\end{multline*}
while for $\beta'\ge \beta$, we have
\[
\int_0^\beta \psi(-t)\log\frac1t\, {\rm d}t \le \psi(0)\,
\bigl( \beta + \beta\log\frac1{\beta} \bigr)
\le \log\frac1{\beta}+1.
\]
That is,
\[
\int_{-\beta}^{\beta'} \log\frac1{|t|}\, \psi (t)\, {\rm d}t
\le 3\log\frac1{\beta\cdot \beta'} + 2\,,
\]
proving~\eqref{eq:3} and completing
the proof of Lemma~\ref{lemma:discretization}. \hfill $\Box$

\section{Riesz products}
Our last results concern with
a family of singular continuous measures introduced by F.~Riesz and
called the Riesz products. These measures  have a variety
of applications in harmonic analysis, see e.g.~\cite[\S13]{Mattila} and the references
therein. Our attention to the Riesz products in the context of this work was attracted by a discussion of Khruschev's work in
\cite[Section~2.11]{Simon}.

To define the Riesz products,
consider a sequence of probability measures
$$
{\rm d} \rho_n(e^{{\rm i}\theta}) = \prod_{j=0}^n\bigl(1+\al_j\cos(\ell_j\theta)\bigr)\,
\frac{{\rm d}\theta}{2\pi},
$$
where $-1\le \al_j \le 1$, and
$\ell_j$ are positive integers such that $\ell_{j+1}\ge 3\ell_j$.
The sequence of measures $\rho_n$ has a weak limit $\rho$ called
{\em the Riesz product}. The measure $\rho$ is singular continuous iff
$$
\sum\limits_{j=0}^\infty \al^2_j=\infty
$$
(otherwise, it is absolutely continuous).

\begin{theorem}\label{thm:RieszProd}
Let $\rho$ be a Riesz product generated by the sequences
$(\alpha_j)$ and $(\ell_j)$, and let $N_n =\sum_{j=0}^n \ell_j$. Then
\[
\prod_{j=0}^n \frac12\, \Bigl( 1+\sqrt{1-\alpha_j^2}\, \Bigr)
\le e_{N_n}(\rho)^2 \le
\prod_{j=0}^n \Bigl( 1 - \frac{\alpha_j^2}4\, \Bigr).
\]
\end{theorem}

In particular, for $\alpha_j\to 0$, we have
\[
2\log e_{N_n}(\rho) = - \frac14 \sum_{j=0}^n \alpha_j^2
+ O\Bigl( \sum_{j=0}^n \alpha_j^4\, \Bigr),
\]
while, for $\alpha_j =1$, $j\in\bZ_+$, we get
\[
-(n+1) \log 2 \le 2 \log e_{N_n}(\rho) \le
-(n+1) \log\frac43\,.
\]

\subsection{Proof of Theorem~\ref{thm:RieszProd}}
First, we note that the moments of the measures $\rho$ and $\rho_n$ coincide up to the order $N_n=\sum_{j=0}^n \ell_j$. So the corresponding orthogonal polynomials (as well as their $L^2(\rho)$- and $L^2(\rho_n)$-norms) coincide too:
$Q_{N_n}(\rho) = Q_{N_n}(\rho_n)$, and $e_{N_n}(\rho)=e_{N_n}(\rho_n)$.

\subsubsection{Proof of the lower bound:}
The proof is straightforward and uses a familiar integral
\[
\int_{-\pi}^{\pi} \log\bigl( 1+\alpha\cos\theta\bigr)\, \frac{{\rm d}\theta}{2\pi}
= \log \Bigl( \frac12 \bigl(1+\sqrt{1-\al^2}\bigr)\, \Bigr).
\]
Since the measure $\rho_n$ has a convergent logarithmic integral,
by Szeg\H{o}'s theorem, for every $k\in\bN$,  we have
\begin{multline*}
\log e_k(\rho_n)\ge \frac 12 \int_{-\pi}^\pi
\log \biggl\{\prod_{j=0}^n\bigl(1+\al_j\cos\(\ell_j\theta\)\bigr)
	\biggr\} \frac{{\rm d}\theta}{2\pi} \\
= \frac 12\sum_{j=0}^n \log\Bigl( \frac12 \bigl(1+\sqrt{1-\al_j^2}\bigr) \Bigr),
\end{multline*}
whence,
$$
\log e_{N_n}(\rho) = \log e_{N_n}(\rho_n)\ge
\frac 12\, \sum_{j=0}^n \log\Bigl( \frac12 \bigl(1+\sqrt{1-\al_j^2}\bigr) \Bigr),
$$
proving the lower bound. \hfill $\Box$

\subsubsection{Proof of the upper bound:}
Consider the monic polynomial
$$P_{N_n}(z):=\prod_{j=0}^n\left(z^{l_j}-\al_j/2\right)$$
of degree $N_n$.
Then
\begin{align*}
e_{N_n}^2(\rho) &= e_{N_n}^2(\rho_{n}) \le \|P_{N_n}\|^2_{L^2(\rho_n)} \\
&=\int_{-\pi}^\pi \prod_{j=0}^n
|e^{{\rm i}\ell_j \theta} -\tfrac12\, \al_j|^2\, (1+\al_j\cos(\ell_{j}\theta))\,
\frac{{\rm d}\theta}{2\pi}
\\
&=\int_{-\pi}^{\pi}  \prod_{j=0}^n
\bigl( 1-\tfrac14\, \al_j^2
+ \tfrac18\, \al_j^3 (e^{{\rm i}\ell_j\theta} + e^{-{\rm i}\ell_j\theta})
- \tfrac14\, \al_j^2 (e^{2{\rm i}\ell_j\theta} + e^{-2{\rm i}\ell_j\theta}) \bigr)\,
\frac{{\rm d}\theta}{2\pi}.
\end{align*}
Observe that due to the growth condition
$ \ell_{j+1}\ge 3\ell_j$, the constant term of the product under the
integral sign, and hence, the whole integral on the RHS is equal to
$$
\prod_{j=0}^n\left(1-\frac {\al_j^2}4\right).
$$
This completes the proof of the upper bound. \hfill $\Box$

\bigskip
\bigskip
\bigskip
\bigskip
\bigskip
\bigskip
\bigskip
\bigskip
\bigskip
\bigskip
\medskip

\noindent{A.B.: Institut de Mathematiques de Marseille,
Aix Marseille Universit\'e, CNRS, Centrale Marseille, I2M, Marseille, France
\newline {\tt alexander.borichev@math.cnrs.fr}
\smallskip
\newline\noindent A.K.: Department of Mathematics and Mechanics,
St. Petersburg State University, St. Petersburg, Russia
\newline {\tt a.kononova@spbu.ru}
\smallskip
\newline\noindent M.S.:
School of Mathematics, Tel Aviv University, Tel Aviv, Israel
\newline {\tt sodin@tauex.tau.ac.il}
}

\end{document}